\documentclass[A4paper,12pt]{article}
\usepackage{amssymb,amsmath}
\usepackage{graphicx}
\newcount\exnom\exnom=0

\long\def\exo/{\vspace{0.05cm} \noindent\advance\exnom by1{\bf
{\the\exnom}.}}

\newcommand{\ds}{\displaystyle}

\begin{document}

\thispagestyle{empty}

\begin{center} \Large
\textbf{JORDAN DERIVATIONS ON SEMIRINGS OF TRIANGULAR MATRICES}
\end{center}

\begin{center} {\bf Dimitrinka Vladeva}\\ University of forestry, bul.$Kl. Ohridski$ 10, Sofia 1000, Bulgaria\\ E-mail: d$\_$vladeva@abv.bg
\end{center}

\begin{quote} \textbf{Abstract}

{\small We explore Jordan derivations of triangular matrices with entries from an additively idempotent semiring. The main result states that for any matrix $A$ over additively idempotent semiring, if we put all the elements of the family of dense submatrices of $A$ to be zeroes, we find a derivative of $A$. The set of derivations of this type is established.}


\end{quote}

\vspace{5mm}

\noindent{\large \bf 1$\;$ Introduction and preliminaries}
\vspace{3mm}

Boolean algebras, fuzzy algebras, bounded distributive lattices, inclines and other algebraic structures are special cases of additively idempotent semirings. These semirings are useful tools in diverse areas such as  automata theory, information systems, dynamic programming and decision theory, see \cite{Golan}. The most interesting additively idempotent semirings are max-plus (or min-plus) algebra applied in modeling network, in language theory, in computer sciences and in idempotent analysis applied in mathematical physics. The techniques of matrices over semirings of this type are well--studied and have important applications, see \cite{But01}, \cite{But02}, \cite{But03} and \cite{But04}.

In 1957 Herstein, see \cite{Her},  proved that arbitrary Jordan derivation from a prime ring of characteristic which is
not 2 into itself is a derivation. This result has been extended by Bre$\check{\mbox{s}}$ar in \cite{Bre01} and \cite{Bre02}.
More recently,  Benkovi$\check{\rm c}$ in \cite{Benk}, studied Jordan derivations in triangular matrices and proved, in particular, that there are no proper Jordan derivations from ring of triangular matrices over commutative ring into itself.

The author studied some derivations in triangular  and other types of matrices over semirings in \cite{Vla01}. Jordan derivations in finite endomorphism semirings are considered in \cite{Vla02}.
Since  considered derivations in this paper are generated by $(0,1)$-- matrices, the reader is referred to \cite{Brualdi}.

In the paper, we deal exclusively with the matrices $E_{ij} = (a_{ij})$, where
$$ a_{ij} =  \left\{ \begin{array}{rl} 1, & \mbox{if}\;\; i = j,\\ 0, & \mbox{if}\;\; i \neq j \end{array} \right. $$ and $0$ and $1$ are zero and identity of some semiring. Evidently
$$E_{ik}\,E_{lj} = \left\{ \begin{array}{rl} E_{ij}, & \mbox{if}\;\; k = l,\\ 0, & \mbox{if}\;\; k \neq l \end{array} \right. .$$

We denote by $UTM_n(S)$ the semiring of upper triangular matrices with entries from some semiring $S$.
\vspace{5mm}

\noindent{\large \bf 2$\;$ Jordan derivations generated by $\mathbf{(0,1)}$ -- matrices}
\vspace{3mm}

Let $A = (a_{ij}) \in UTM_n(S_0)$, where $S_0$ is an additively idempotent semiring, that is $\displaystyle A = \sum_{i \leq j} a_{ij}E_{ij}$, where $a_{ij} \in S_0$, $i,j = 1, \ldots n$.
Let $\overline{D}_{k} = E_{11} + \cdots + E_{kk}$, where $1 \leq k \leq n$. We obtain
$$A\,\overline{D}_{k} = a_{11}E_{11} + a_{12}E_{12} + \cdots + a_{1k}E_{1k} + a_{22}E_{22} + \cdots + a_{2k}E_{2k} + \cdots + a_{kk}E_{kk}.$$

 Similarly we find
$$\overline{D}_{k}\,A = a_{11}E_{11} + a_{12}E_{12} + \cdots + a_{1n}E_{1n}  + \cdots + a_{kk}E_{kk} + \cdots + a_{kn}E_{kn}.$$

Since $A\,\overline{D}_{k}$ is a submatrix of $\overline{D}_{k}\,A$ for Jordan product, it follows
$$A\circ \overline{D}_{k} = A\,\overline{D}_{k} + \overline{D}_{k}\,A = \overline{D}_{k}\,A, \; \mbox{or}$$
$$A\circ \overline{D}_{k} = \left(
                           \begin{array}{ccccccc}
                             a_{11} & a_{12} & \cdots & a_{1k} &  \cdots & a_{1n} \\
                             0 & a_{22} & \cdots & a_{2k} &  \cdots & a_{2n} \\
                             \cdots & \cdots & \ddots & \cdots &  \cdots & \cdots \\
                             0 & 0 & \cdots & a_{kk}\; &  \cdots & a_{kn}\\
                             0 & 0 & \cdots & 0 &  \cdots & 0 \\
                             \cdots & \cdots & \cdots &  \cdots & \cdots & \cdots \\
                             0 & 0 & \cdots & 0 &  \cdots & 0 \\
                           \end{array}
                         \right).$$

Let $\displaystyle B = \sum_{i \leq j} b_{ij}E_{ij} \in UTM_n(S_0)$. Then, it follows
$$ B\circ \overline{D}_{k} = \overline{D}_{k}\,B  \;\mbox{and}\; (A\,B)\circ \overline{D}_{k} = \overline{D}_{k}\,(A\,B).$$

Now we obtain
$$(A\circ \overline{D}_{k})B + A(B\circ \overline{D}_{k}) = \overline{D}_{k}\,A\,B + A\,\overline{D}_{k}\,B = $$ $$= (\overline{D}_{k}\,A + A\,\overline{D}_{k})B = \overline{D}_{k}\,(A\,B) = (A\,B)\circ \overline{D}_{k}.$$

It is easy to see that $(A + B)\circ \overline{D}_{k} = \overline{D}_{k}\,(A + B) = A\circ \overline{D}_{k} + B\circ \overline{D}_{k}$.
Thus we prove

\vspace{3mm}

\textbf{Proposition 1.} \textsl{The map  $\delta_k : UTM_n(S_0) \rightarrow UTM_n(S_0)$, where $S_0$ is an additively idempotent semiring and  $\delta_k(A) = A\circ \overline{D}_{k}$ for any matrix $A \in UTM_n(S_0)$, is a derivation}.

 \vspace{1mm}

 It is well known that for $k = n$ the map $\delta_n = i$ is an identity map, which is a derivation.

 \vspace{1mm}

 Since for any matrix $A$, it follows
 $$\delta_k(A) + \delta_{\ell}(A) = \delta_{\ell}(A), \; \mbox{and} \; \delta_\ell\delta_k(A) = \delta_k(\delta_\ell(A)) = \delta_k(A),$$
 where $k \leq \ell$, we can write
 $$\delta_k + \delta_\ell = \delta_\ell + \delta_k = \delta_\ell,$$
 $$\delta_k\, \delta_\ell = \delta_\ell\, \delta_k = \delta_k$$
 for $k \leq \ell$.

 Let us denote by $\overline{\mathcal{D}}$ the set of derivations $\delta_k$, where $k = 1, \ldots, n$. Then $(\overline{\mathcal{D}},+,.)$ is a semiring with a zero and smallest element $\delta_1$ and identity and greatest element the identity map $\delta_n$.

 Thus $(\overline{\mathcal{D}},+,.)$ is an additively idempotent and also a multiplicatively idempotent semiring.

 \vspace{5mm}

 Let $\underline{D}_{\,m} = E_{n-m+1\,n-m+1} + \cdots + E_{nn}$, where $1 \leq m \leq n$. We find
$$A\,\underline{D}_{\,m} = a_{1\,n-m+1}E_{1\,n-m+1}  + \cdots + a_{1n}E_{1n}  + \cdots + $$ $$+ a_{n-m+1\,n-m+1}E_{n-m+1\,n-m+1} + \cdots + a_{n-m+1\,n}E_{n-m+1\,n} +$$ $$+ a_{n-m\,n-m}E_{n-m\,n-m} + \cdots + a_{n-m\,n}E_{n-m\,n} + \cdots + a_{nn}E_{nn}.$$

 Similarly we obtain
 $$\underline{D}_{\,m}\,A = a_{n-m+1\,n-m+1}E_{n-m+1\,n-m+1} + \cdots + a_{n-m+1\,n}E_{n-m+1\,n} +$$
 $$+ a_{n-m\,n-m}E_{n-m\,n-m} + \cdots + a_{n-m\,n}E_{n-m\,n} + \cdots + a_{nn}E_{nn}.$$

Since $\underline{D}_{\,m}\,A$ is a submatrix of $A\,\underline{D}_{\,m}$ for Jordan product, it follows
$$A\circ \underline{D}_{\,m} = A\,\underline{D}_{\,m} + \underline{D}_{\,m}\,A = A\,\underline{D}_{\,m}, \; \mbox{or}$$
$$A\circ \underline{D}_{\,m} = \left(
                            \begin{array}{ccclcl}
                              0 & \cdots & 0 & a_{1\,n-m+1} & \cdots & a_{1n} \\
                              0 & \cdots & 0 & a_{2\,n-m+1} & \cdots & a_{2n} \\
                              \vdots & \vdots & \vdots & \vdots & \vdots & \vdots \\ \vspace{2mm}
                              0 & \cdots & 0 & a_{n-m+1\,n-m+1} & \cdots & a_{n-m+1\,n} \\ \vspace{1mm}
                              \cdots & \cdots & \cdots & \cdots & \ddots & \cdots \\
                              0 & \cdots & 0 & 0 & \cdots & a_{nn} \\
                            \end{array}
                          \right).$$

Let $\displaystyle B \in UTM_n(S_0)$. Then, it follows
$$ B\circ \underline{D}_{\,m} = B\,\underline{D}_{\,m} \;\mbox{and}\; (A\,B)\circ \underline{D}_{\,m} = (A\,B)\,\underline{D}_{\,m}.$$

Now we find
$$(A\circ \underline{D}_{\,m})B + A(B\circ \underline{D}_{\,m}) = A\,\underline{D}_{\,m}\,B + A\,B\,\underline{D}_{\,m} = $$ $$= A\,(\underline{D}_{\,m}\,B + B\,\underline{D}_{\,m}) = (A\,B)\underline{D}_{\,m} = (A\,B)\circ \underline{D}_{\,m}.$$

It is easy to see that $(A + B)\circ \underline{D}_{\,m} = (A + B)\,\underline{D}_{\,m} = A\circ \underline{D}_{\,m} + B\circ \underline{D}_{\,m}$. Thus we prove

\vspace{3mm}

\textbf{Proposition 2.} \textsl{The map  $d_m : UTM_n(S_0) \rightarrow UTM_n(S_0)$, where $S_0$ is an additively idempotent semiring and  $d_m(A) = A\circ \underline{D}_{\,m}$ for any matrix $A \in UTM_n(S_0)$, is a derivation}.

 \vspace{2mm}

  For $m = n$ the map $d_n = i$ is an identity map, which is a derivation.

 \vspace{2mm}

 Since for any matrix $A$, it follows
 $$d_\ell(A) + d_{m}(A) = d_{m}(A), \; \mbox{and} \; d_\ell d_m(A) = d_m(d_\ell(A)) = d_{\ell}(A),$$
 where $\ell \leq m$, we can write
 $$d_m + d_\ell = d_\ell + d_m = d_m,$$
 $$d_m d_\ell = d_\ell d_m = d_\ell$$
 for $\ell \leq m$.

 Let us denote by $\underline{\mathcal{D}}$ the set of derivations $d_m$, where $m = 1, \ldots, n$. Then $(\underline{\mathcal{D}},+,.)$ is a semiring with a zero and smallest element $d_n$ and identity and greatest element the identity map $d_1$.

 Hence,  $(\underline{\mathcal{D}},+,.)$ is an additively idempotent and also a multiplicatively idempotent semiring.

Let ${\mathcal{D}} = \overline{\mathcal{D}} \cup \underline{\mathcal{D}}$. Since elements of semirings $(\overline{\mathcal{D}},+,.)$ and $(\underline{\mathcal{D}},+,.)$ are derivations in the semiring $UTM_n(S_0)$ we can add them and their sums are derivations. Thus $\delta_k + d_m \in \mathcal{D}$  are derivations for any $1 \leq k \leq n$ and $1 \leq m \leq n$.
For matrix $A = (a_{ij}) \in UTM_n(S_0)$ and $k + m < n$, it follows
$$(\delta_k + d_m)(A) = \left(
                            \begin{array}{cccccclcl}
                              a_{11} & \cdots & a_{1k} & a_{1\,k+1} & \cdots & a_{1\,n-m} & a_{1\,n-m+1} & \cdots & a_{1n} \\
                              \vdots & \ddots & \vdots & \vdots & \vdots & \vdots & \vdots & \vdots & \vdots\\
                              0 & \cdots & a_{kk} & a_{k\,k+1} & \cdots & a_{k\,n-m} & a_{k\,n-m+1} & \cdots & a_{k\,n} \\
                              0 & \cdots & 0 & \mathbf{0} & \cdots & \mathbf{0} & a_{k+1\,n-m+1} & \cdots & a_{k+1\,n} \\
                               \vdots & \vdots & \vdots & \vdots & \ddots & \vdots & \vdots & \vdots & \vdots\\
                              0 & \cdots & 0 & \mathbf{0} & \cdots & \mathbf{0} & a_{n-m\,n-m+1} & \cdots & a_{n-m\,n} \\ \vspace{2mm}
                              0 & \cdots & 0 & {0} & \cdots & {0} & a_{n-m+1\,n-m+1} & \cdots & a_{n-m+1\,n} \\
                              \cdots & \cdots & \cdots & \cdots & \cdots & \cdots &\cdots & \ddots & \cdots \\
                              0 & \cdots & 0 & 0 & \cdots & 0 & 0& \cdots & a_{nn} \\
                            \end{array}
                          \right).$$

   The semigroup $(\mathcal{D},+)$ has no zero element. We can define a particular order in $(\mathcal{D},+)$ by rules
 $$\delta_\ell \leq \delta_k + d_m \; \mbox{if}\; \ell \leq k,$$
 $$\;d_\ell \leq \delta_k + d_m \; \mbox{if}\; \ell \leq m,$$
 $$\delta_{k_1} + d_{m_1}  \leq \delta_k + d_m \; \mbox{if}\; k_1 \leq k \; \mbox{and}\; m_1 \leq m.$$

A principal dense submatrix of a matrix $A = (a_{ij}) \in  UTM_n(S_0)$ is a square submatrix of $A$ in which the main diagonal consists of the elements $a_{ii}, \ldots, a_{jj}$, where $i, \ldots, j$ are consecutive numbers and $1 \leq i \leq j \leq n$. The matrix $(\delta_k + d_m)(A)$ from above has a principal dense submatrix with main diagonal $0_{k+1\,k+1}, \ldots, 0_{n-m\,n-m}$ and all entries equal to zero.

\vspace{3mm}

\textbf{Theorem 1.} \textsl{ Let for any matrix $A = (a_{ij}) \in  UTM_n(S_0)$, where $S_0$ is an additively idempotent semiring, the matrix $\delta(A)   \in  UTM_n(S_0)$ have the same elements $a_{ij}$, $1 \leq i \leq n$, $1 \leq j \leq n$ except the elements of a principal dense submatrix of $A$, which are zeroes. Then the map  $\delta : UTM_n(S_0) \rightarrow UTM_n(S_0)$ is a derivation. The number of all these derivations $\delta$ is $\ds \binom{n+1}{2}$.}

\emph{Proof.} The statement that $\delta$ is a derivations follows from Proposition 1, Proposition 2 and constructions of derivations $\delta_k + d_m$, where $1 \leq k \leq n$ and $1 \leq m \leq n$.

Let $s$ be the order of the principal dense submatrix of $A$ with zero elements.

When $s = 0$ the derivation $\delta$ is the identity map.

When $s = 1$ there are $n$ derivations (namely $d_{n-1}, \delta_1 + d_{n-2}, \ldots, \delta_{n-2} + d_1, d_{n-1}$).

When $s = 2$ we obtain $n - 1$ derivations (namely $d_{n-2}, \delta_1 + d_{n-3}, \ldots, \delta_{n-3} + d_1, \delta_{n-2}$).

For arbitrary $s = i$, where $1 \leq i \leq n-1$,  there are $n - i + 1$ derivations.

When $s = n-1$ we find two derivations ($\delta_1$ and $d_1$).

Thus the number of all derivations $\delta$ is $\ds 1 + n + n - 1 + \cdots + 2 = \binom{n+1}{2}$.

\vspace{5mm}

\noindent{\large \bf 3$\;$ Products of derivations}

\vspace{3mm}

The product $\delta_k d_m$, where $1 \leq k \leq n$ and $1 \leq m \leq n$ is well defined by rule $\delta_k d_m(A) = d_m(\delta_k(A))$ for any $A  \in UTM_n(S_0)$, but in the general case it is not  a derivation as we see in the following example.

\vspace{2mm}

 \textbf{\emph{Example}} Let us consider the matrix $A = \left(
                                        \begin{array}{ccc}
                                          a_{11} & a_{12} & a_{13} \\
                                          0 & a_{22} & a_{23} \\
                                          0 & 0 & a_{33} \\
                                        \end{array}
                                      \right)$ with entries ${a_{ij} \in S_0}$, where $S_0$ is an additively idempotent semiring.
\vspace{3mm}

                                      Then we obtain $\delta_1(A) = \left(
                                        \begin{array}{ccc}
                                          a_{11} & a_{12} & a_{13} \\
                                          0 & 0 & 0 \\
                                          0 & 0 & 0 \\
                                        \end{array}
                                      \right)$ and $d_3(A) = \left(
                                        \begin{array}{ccc}
                                          0 & 0 & a_{13} \\
                                          0 & 0 & a_{23} \\
                                          0 & 0 & a_{33} \\
                                        \end{array}
                                      \right)$.

                                       Thus we find $\delta_1 d_3 (A) = d_3(\left(
                                        \begin{array}{ccc}
                                          a_{11} & a_{12} & a_{13} \\
                                          0 & 0 & 0 \\
                                          0 & 0 & 0 \\
                                        \end{array}\right)) = a_{13}E_{13}$.
\vspace{3mm}

                                         For any matrix $B = \left(
                                        \begin{array}{ccc}
                                          b_{11} & b_{12} & b_{13} \\
                                          0 & b_{22} & b_{23} \\
                                          0 & 0 & b_{33} \\
                                        \end{array}
                                      \right)$ with entries $b_{ij} \in S_0$, we obtain that $\delta_1 d_3 (B) = b_{13}E_{13}$ and also ${\delta_1 d_3 (AB)\vphantom{\int^A} = a_{11}b_{13} + a_{12}b_{23} + a_{13}b_{33}}$.
Now we calculate
$$\delta_1 d_3 (A)\,B + A\,\delta_1 d_3 (B) = a_{13}E_{13}B + A\,b_{13}E_{13} = (a_{13}b_{33} + a_{11}b_{13})E_{13} \neq \delta_1 d_3 (AB).$$

\vspace{1mm}

\textbf{Theorem 2.} \textsl{Let $\delta_k, d_m \in \mathcal{D}$, where $1 \leq k \leq n$ and $1 \leq m \leq n$.
The map $\delta_k d_m = d_m\delta_k$ is a derivation if and only if  $\delta_k + d_m$ is the identity map.}

\emph{Proof.} Let $k + m \geq n \Leftrightarrow \overline{D}_{k} + \underline{D}_{m} = E \Leftrightarrow \delta_k + d_m = i$, where $E$ and $i$ are the identity matrix and the identity map.

In propositions 1 and 2 we prove that $\delta_k(A) = \overline{D}_{k}A$ and $d_m(A) = A\underline{D}_{m}$. Then, it follows $\delta_k d_m(A) = d_m(\delta_k(A)) = \overline{D}_{k}A\underline{D}_{m}$. Hence we obtain
$$(\delta_k d_m(A))B = \overline{D}_{k}A\underline{D}_{m}B, \; A(\delta_k d_m(B)) = A\overline{D}_{k}B\underline{D}_{m} \; \mbox{and}\; \delta_k d_m(AB) = \overline{D}_{k}AB\underline{D}_{m}.$$

The last matrix has the following representation:
$$\overline{D}_{k}AB\underline{D}_{m} = \overline{D}_{k}A(\overline{D}_{k} + \underline{D}_{m})B\underline{D}_{m} = \overline{D}_{k}A\overline{D}_{k}B\underline{D}_{m} + \overline{D}_{k}A\underline{D}_{m}B\underline{D}_{m}.$$

But $\overline{D}_{k}A\overline{D}_{k} = (a_{11}E_{11} + a_{12}E_{12} + \cdots + a_{1n}E_{1n}  + \cdots + a_{kk}E_{kk} + \cdots + a_{kn}E_{kn})(E_{11} + \cdots + E_{kk}) = a_{11}E_{11} + a_{12}E_{12} + \cdots + a_{1k}E_{1k} +  \cdots + a_{kk}E_{kk} = A\overline{D}_{k}$. Similarly we find $\underline{D}_{m}B\underline{D}_{m} = \underline{D}_{m}B.$

Hence $\overline{D}_{k}AB\underline{D}_{m} = A\overline{D}_{k}B\underline{D}_{m} + \overline{D}_{k}A\underline{D}_{m}B$, that is $\delta_k d_m(AB) =  A(\delta_k d_m(B)) + (\delta_k d_m(A))B$. Since the map $\delta_k d_m$ is evidently linear, it follows that $\delta_k d_m$ is a derivation.

Let $k + m < n$, i.e. the derivation $\delta_k + d_m$ is not the identity map. Let us compare the elements in the first row and $n$-th column of the matrices $(\delta_k d_m(A))B + A(\delta_k d_m(B))$ and $\delta_k d_m(AB)$.

Since $(\delta_k d_m(A))B = (\overline{D}_{k}A)(\underline{D}_{m}B) = (a_{11}E_{11} + a_{12}E_{12} + \cdots + a_{1n}E_{1n}  + \cdots + a_{kk}E_{kk} + \cdots + a_{kn}E_{kn})(b_{n-m+1\,n-m+1}E_{n-m+1\,n-m+1} + \cdots + b_{n-m+1\,n}E_{n-m+1\,n} +
  b_{n-m\,n-m}E_{n-m\,n-m} + \cdots + b_{n-m\,n}E_{n-m\,n} + \cdots + b_{nn}E_{nn})$, it follows that the element in the first row and $n$-th column of this matrix is $\ds \sum_{i = n-m+1}^n a_{1i}b_{in}$.

  Analogously $A(\delta_k d_m(B)) = (A\overline{D}_{k})(B\underline{D}_{m}) = (a_{11}E_{11}  + \cdots + a_{1k}E_{1k} + a_{22}E_{22} + \cdots + a_{2k}E_{2k} + \cdots + a_{kk}E_{kk})(b_{1\,n-m+1}E_{1\,n-m+1}  + \cdots + b_{1n}E_{1n}  + \cdots +  b_{n-m+1\,n-m+1}E_{n-m+1\,n-m+1} + \cdots + b_{n-m+1\,n}E_{n-m+1\,n} + b_{n-m\,n-m}E_{n-m\,n-m} + \cdots + b_{n-m\,n}E_{n-m\,n} + \cdots + b_{nn}E_{nn})$ and then the element in the first row  and $n$-th column of this matrix is $\ds \sum_{i = 1}^k a_{1i}b_{in}$.

  The element in the first row  and $n$-th column of  matrix $\delta_k d_m(AB)$ is $\ds \sum_{i = 1}^n  a_{1i}b_{in}$.
  However in general case $\ds \sum_{i = n-m+1}^n a_{1i}b_{in} + \sum_{i = 1}^k a_{1i}b_{in} \neq \sum_{i = 1}^n  a_{1i}b_{in}$, because at least  the summand $a_{1\,n-m}b_{n-m\,n}$ not occur in the left sum.

\vspace{2mm}

  Immediately from the last theorem follows

\vspace{3mm}

\textbf{Corollary 1.} \textsl{For for any matrix $A = (a_{ij}) \in  UTM_n(S_0)$, where $S_0$ is an additively idempotent semiring and derivation $\delta_k \in \mathcal{D}$, where $k = 1, \ldots, n-1$, the product $\delta_kd_{n-k}$ is a derivation and}
$$\delta_kd_{n-k}(A) = \left(
                            \begin{array}{cccccc}
                              \mathbf{0} & \cdots & \mathbf{0} & a_{1\,k+1} & \cdots & a_{1\,n}\\
                              \vdots & \ddots & \vdots & \vdots & \vdots & \vdots\\
                              \mathbf{0} & \cdots & \mathbf{0} & a_{k\,k+1} & \cdots & a_{k\,n}\\
                                  0 & \cdots & 0 & \mathbf{0} &  \cdots &\mathbf{0}\\
                              \vdots & \vdots & \vdots & \vdots & \ddots & \vdots\\
                              0 & \cdots & 0 & \mathbf{0} & \cdots & \mathbf{0}\\
                            \end{array}
                          \right).$$

Derivations $\delta_k$ and $d_{n-k}$ are called complementary derivations. In [12] we prove that the map $D : UTM_n(S_0) \rightarrow UTM_n(S_0)$, where $S_0$ is an additively idempotent semiring, such that
$$D(A) = A\backslash\mbox{diag}(a_{11}, \ldots, a_{nn})$$
for any matrix $A = (a_{ij}) \in  UTM_n(S_0)$ is a derivation. The derivation $D$ can be represented as a sum of products of all complementary derivations $\delta_k$ and $d_{n-k}$:
$$D = \delta_1d_{n-1} + \delta_2d_{n-2} + \cdots + \delta_{n-1}d_1.$$

 Immediately from the last theorem follows

\vspace{3mm}

\textbf{Corollary 2.} \textsl{For  any matrix $A = (a_{ij}) \in  UTM_n(S_0)$, where $S_0$ is an additively idempotent semiring and derivation $\delta_k \in \mathcal{D}$, where $k = 1, \ldots, n-1$, the product $\delta_kd_{n-k + p}$, where $1 \leq p \leq k$, is a derivation and}
$$\delta_kd_{n-k + p}(A) = \left(
                            \begin{array}{ccclcllcl}
                              \mathbf{0} & \cdots & \mathbf{0} & a_{1\,k-p+1} &\cdots & a_{1\,k} & a_{1\,k+1} & \cdots & a_{1\,n}\\
                              \vdots & \ddots & \vdots & \vdots & \vdots & \vdots  & \vdots  & \vdots  & \vdots\\
                              \mathbf{0} & \cdots & \mathbf{0} & a_{k-p\,k-p+1} &\cdots & a_{k-p+1\,k} & a_{k-p+1\,k+1} & \cdots & a_{k-p+1\,n}\\
                               {0} & \cdots & {0} & a_{k-p+1\,k-p+1} &\cdots & a_{k-p+1\,k} & a_{k-p+1\,k+1} & \cdots & a_{k-p+1\,n}\\
                              \vdots & \vdots & \vdots & \vdots & \ddots \;\;\;\;\; & \vdots  & \vdots  & \vdots  & \vdots\\
                              0 & \cdots & 0 & 0 & \cdots & a_{kk} & \;\;\;a_{k\,k+1} & \cdots & a_{kn}\\
                              0 & \cdots & 0  &    0 & \cdots & 0 & \;\;\;\;\;\mathbf{0} &  \cdots &\mathbf{0}\\
                             \vdots & \vdots & \vdots & \vdots & \vdots & \vdots & \vdots & \ddots \;\;\;\;\; & \vdots\\
                             0 & \cdots & 0  & 0 & \cdots & 0 & \;\;\;\;\; \mathbf{0} & \cdots & \mathbf{0}\\
                            \end{array}
                          \right).$$

Let us consider the matrix $A = (a_{ij}) \in  UTM_n(S_0)$, where $S_0$ is an additively idempotent semiring. Denote by $A(i_k,n_k)$ the principal dense submatrix of $A$ with main diagonal
$$a_{i_k\,i_k} \; \cdots \; a_{i_k + n_k - 1\,i_k + n_k - 1}.$$

An arbitrary finite set of principal dense submatrices $A(i_k,n_k)$, where $k = 1, \ldots, s$, without common elements is called a family.

\vspace{3mm}

\textbf{Theorem 3.} \textsl{ Let for any matrix $A = (a_{ij}) \in  UTM_n(S_0)$, where $S_0$ is an additively idempotent semiring, the matrix $\delta(A)   \in  UTM_n(S_0)$ have the same elements $a_{ij}$, $1 \leq i \leq n$, $1 \leq j \leq n$ except the elements of a family of principal dense submatrices $A(i_k,n_k)$, where $k = 1, \ldots, s$, which are zeroes. Then the map  $\delta : UTM_n(S_0) \rightarrow UTM_n(S_0)$ is a derivation. The number of all these derivations $\delta$ is $\ds 2^n$.}

\vspace{1mm}

\emph{Proof.} Let on the main diagonal of the matrix $\delta(A)$ there are only zeroes. In other words the family of principal dense submatrices
$$A(1,n_1), A(n_1 + 1,n_2), \ldots, A(n_s + 1,n)$$
covers the main diagonal. Then, using the complemented derivations (see Corrolary 1) we obtain
$$\delta = \delta_{n_1}d_{n-n_1} + \delta_{n_2}d_{n-n_2} + \cdots + \delta_{n_s}d_{n-n_s}.$$

Let the nonzero elements on the main diagonal of the matrix $\delta(A)$ are:
$$\begin{array}{lcl} a_{11}, & \cdots, & a_{i_1-1\,i_1-1},\\
a_{i_1 + n_1\,i_1 + n_1}, & \cdots, & a_{i_2-1\,i_2-1},\\
a_{i_2 + n_2\,i_2 + n_2}, & \cdots, & a_{i_3-1\,i_3-1},\\
\cdots & \cdots & \cdots\\
a_{i_{s-1} + n_{s-1}\,i_{s-1} + n_{s-1}}, & \cdots, & a_{i_{s}-1\,i_{s}-1}\\
\end{array}.\eqno{(1)}$$

This construction of the set of nonzero elements corresponds to the family of principal dense submatrices
$$A(i_1,n_1), A(i_2,n_2), \ldots, A(i_{s-1},n_{s-1}), A(i_{s},n_{s}), \; \mbox{where}\; i_{s} + n_{s} - 1 = n$$
with zero elements. Then, using the  derivations $\delta_kd_{n-k + p}$ from Corrolary 2 we find
$$\delta = \delta_{i_1 - 1} + \delta_{i_2 - 1}d_{n - i_1 - n_1 + 1} + \cdots + \delta_{i_{s} - 1}d_{n - i_{s-1} - n_{s-1} + 1}. \eqno{(2)}$$

By the similar reasonings we consider the case when the first diagonal element of the matrix $\delta(A)$ is zero. So, when nonzero elements of the main diagonal are those in (1) without elements on first row, from (2) we obtain
$$\delta =  \delta_{i_2 - 1}d_{n - i_1 - n_1 + 1} + \cdots + \delta_{i_{s} - 1}d_{n - i_{s-1} - n_{s-1} + 1}. \eqno{(3)}$$

Analogously, when the last nonzero elements on the main diagonal  of the matrix $\delta(A)$ are
$$a_{i_{s} + n_{s}\,i_{s} + n_{s}}, \cdots, a_{nn}$$
from (2) (or from (3))
we have
$$\delta = \delta_{i_1 - 1} + \delta_{i_2 - 1}d_{n - i_1 - n_1 + 1} + \cdots + \delta_{i_{s} - 1}d_{n - i_{s-1} - n_{s-1} + 1} + d_{n - i_{s} - n_s + 1}$$
or respectively
$$\delta =  \delta_{i_2 - 1}d_{n - i_1 - n_1 + 1} + \cdots + \delta_{i_{s} - 1}d_{n - i_{s-1} - n_{s-1} + 1} + d_{n - i_{s} - n_s + 1}.$$

\noindent{\large \bf 4$\;$ Other derivations}
\vspace{2mm}

In order to obtain derivations different from the above considered we explore the
 max-plus semiring $(\mathbb{R}_{\max},\oplus,\odot)$, where
$\mathbb{R}_{\max} = \mathbb{R} \cup \{- \infty\}$,  $\mathbb{R}$ is the field of real numbers and for any $a, b \in \mathbb{R}$
 $$a \oplus b = \max\{a,b\}\;\;\; \mbox{and}\;\;\; a \odot b = a + b.$$

For any $x \in \mathbb{R}_{\max}$ we define a map $\delta_x : \mathbb{R}_{\max} \rightarrow \mathbb{R}_{\max}$ such that
  $$\delta_x(a) = a\odot x = a + x,$$ where $a \in \mathbb{R}_{\max}$. When $x = 0$ obviously $\delta_x$ is an identity map in $\mathbb{R}_{\max}$.

\vspace{2mm}

\textbf{Proposition 3.}  \textsl{The map $\delta_x : \mathbb{R}_{\max} \rightarrow \mathbb{R}_{\max}$ is a derivation for any $x \in \mathbb{R}_{\max}$.}

\emph{Proof.} We obtain
$$\delta_x(a\oplus b) = (a\oplus b)\odot x = (a\odot x)\oplus(b\odot x) = \delta_x(a)\oplus \delta_x(b).$$

Since $\delta_x(a\odot b) = (a\odot b)\odot x = a + b + x$, $\delta_x(a)\odot b = (a\odot x)\odot b = a + x + b = a + b + x$ and
$a\odot \delta_x(b) = a\odot (b\odot x) = a + b + x$, it follows
$$(\delta_x(a)\odot b)\oplus (a\odot \delta_x(b)) = \max\{(a + b + x),(a + b + x)\} = a + b + x = \delta_x(a\odot b).$$

Let $x, y \in \mathbb{R}_{\max}$. Since $$\delta_y(\delta_x(a)) = \delta_y(a\odot x) = (a\odot x)\odot y = a + x + y = \delta_{x\odot y}(a)$$ for any $a \in \mathbb{R}_{\max}$, it follows $\delta_x . \delta_y = \delta_{x\odot y}(a)$.

Hence for any real number $x$ we obtain that $\delta_x . \delta_{-x}$ is the identity map. It is easy to see that $\delta_x^2 = \delta_x$ for any $x$.
\vspace{2mm}

Let $DER_{\max}$ be the set of all derivations  $\delta_x$ where $x \in \mathbb{R}_{\max}$. Thus we obtain

\vspace{2mm}

\textbf{Corrolary 3.}  \textsl{The set $DER_{\max}$ is an idempotent Abelian group which identity is the identity map.}

It is easy to prove the following fact.

\vspace{2mm}

\textbf{Proposition 4.}  \textsl{Let $(S,+,.)$ be an additively idempotent semiring and maps ${\delta_1 : S \rightarrow S}$ and $\delta_2 : S \rightarrow S$ are derivations. Then $\delta_1 + \delta_2$ defined by $(\delta_1 + \delta_2)(x) = \delta_1(x) + \delta_2(x)$ for any $x \in S$ is a derivation in $S$.}

Immediately follows

\vspace{1mm}

\textbf{Corrolary 4.}  \textsl{$\left(DER_{\max},\oplus,\odot\right)$ is an additively idempotent semifield and all its element are multiplicatively idempotent.}

Now from Proposition 1, [12], for any $\delta_x \in DER_{\max}$ we define the hereditary derivation $\delta^h_x$ in semiring $UTM_n(\mathbb{R}_{\max})$.
In the same way we define derivations in semiring  $UTM_n(S)$, where $S$ is a finite endomorphism semiring. Hence, in general, in $UTM_n(S_0)$, where $S_0$ is an additively idempotent semiring, there are derivations different from considered Jordan derivations.

\end{document}